\documentclass{article}

\usepackage{amsmath}  
\usepackage{amssymb}  
\usepackage{amsfonts} 
\usepackage{epsfig}
\usepackage{url}

\title{ Exponential relaxation data analysis by parametrized regularization of severely ill-posed Fredholm integral equations of the first kind}
\author{V. V. Kryzhniy}
\date{Seattle, September 4, 2023} 

\begin{document}
\maketitle

\begin{abstract}

This paper presents a novel approach to construct regularizing operators for severely ill-posed Fredholm integral equations of the first kind by introducing parametrized discretization. The optimal values of discretization and regularization parameters are computed simultaneously by solving a minimization problem formulated based on a regularization parameter search criterion. The effectiveness of the proposed approach is demonstrated through examples of noisy Laplace transform inversions and the deconvolution of nuclear magnetic resonance relaxation data.
\end{abstract}

\section{Introduction}
The analysis of exponential relaxation data poses a significant challenge in various fields such as experimental physics, chemistry, electrochemistry, and biophysics \cite{Istratov}. This problem involves determining the distribution function $f(t)$ from experimentally measured function $g(s)$ by performing the inversion of the Laplace transformation (\ref{1}):

\begin{equation}
\label{1}
g(s) = \int_0^\infty{\mathrm{e}^{-st}f(t)\mathrm{d}t},
\end{equation}

A similar problem arises in nuclear magnetic resonance relaxometry (NMR) \cite{Kroeker}, where the desired distribution is obtained by deconvolving the following integral equation:
\footnote{Although the equation (\ref{100}) can be written as a Laplace transform with the help of substitution $\tau = t^{-1}$, a direct regularized solution of equation (\ref{100}) is advantageous when $g(s)$ is noisy.}

\begin{equation}
\label{100}
g(s) =\int_0^\infty{\mathrm{e}^{-s/t}f(t)\mathrm{d}t} .
\end{equation}

These equations are inherently ill-posed problems and are extremely sensitive to small perturbations in the right-hand side $g(s)$ \cite{Hansen}. Due to smoothing property of exponential kernel the equations (\ref{1}), (\ref{100}) are severely ill-posed problems \cite{Hansen}.

Following \cite{Varah}, let us illustrate why  these equations are severely ill-posed.
By discretizing the integral equation  (\ref{1}) or (\ref{100}) using an appropriate quadrature formula, we obtain the corresponding matrix equation:

\begin{equation}
\label{3}
K f = g,
\end{equation}

where $K$ is an $m \times n$ matrix, $f$ and $g$ are column vectors of length $n$ and $m$, respectively.

By representing matrix $K$ in equation (\ref{3}) using singular value decomposition (SVD):

\begin{displaymath}
K = U D V^T,
\end{displaymath}

where $U$ and $V$ are orthogonal matrices, $D = \mathrm{diag}(s_i)$, and $s_1 \geq s_2 \geq \dots \geq 0$, we obtain the formal solution of Eq. (\ref{3}):

\begin{equation}
\label{4}
f = K^{-1} g = \sum_{i=1}^n {s_i^{-1}(u_i^T g) v_i}{\mathstrut}.
\end{equation}

Due to division by decreasing singular values $s_i$ in equation (\ref{4}), the formal solution is highly oscillating and unreasonable. Truncated SVD solution regularizes the problem by limiting the number of terms in equation (\ref{4}) by a certain number $n_0$:

\begin{equation}
\label{5}
f_r = K_r^{-1} g = \sum_{i=1}^{n_0} {s_i^{-1}(u_i^T g) v_i} {\mathstrut},
\end{equation}

where $f_r$ and $K_r^{-1}$ represent the regularized solution and regularizing operator, respectively.

The number of terms in equation (\ref{5}) depends on the noise level in the data $\epsilon$ and the rate of decrease of singular values $s_i$ \cite{Hansen, Varah}. For an exponential kernel, the singular values decrease rapidly, resulting in only a few terms in the sum of equation (\ref{5}). Consequently, the regularized solution $f_r$ is represented as a linear combination of a small number of vectors $v_i$:

\begin{equation}
\label{6}
f_r = \frac{\beta_1}{s_1} v_1 + \frac{\beta_2}{s_2} v_2 + \dots,
\end{equation}

where $\beta = U^T g$.

It is evident that such a computed regularized solution is generally inaccurate. This conclusion holds true for all known regularization techniques \cite{Varah}.

\section{Towards regularization of severely ill-posed problems}

All terms in equation (\ref{6}) depend on the SVD of the kernel matrix $K$, which, in turn, depends on the quadrature formula, selected nodes, and the number of points. Although we cannot change the fact that the singular values of the exponential kernel rapidly tend to zero, we can anticipate that a few terms in equation (\ref{6}) will yield more accurate results with carefully tailored discretization.

Thus, for obtaining a reasonable solution of matrix equation (\ref{3}), we need to find an appropriate discretization and value of the regularization parameter.

This idea of introducing flexible discretization seems very natural. The most robust quadrature programs were developed a long time ago \cite{quadpack}, and they are based on fine-tuned partitioning of the interval of integration. Nevertheless, seemingly all modern recommendations for discretizing the integral equations  are based on the standard quadrature formulae, and regularization methods do not include discretization into consideration \cite{Hansen}. That is, currently, discretization and regularization are two disjointed steps for solving integral equations.

The reason for disjointed consideration of discretization and regularization may be attributed to the absence of theory for constructing regularizing operators that include additional parameters along with the regularization parameter \cite{Tikhonov}. At present time, we can rely on a precedent.

For inverting real Laplace transforms, the author has derived an integral form of the regularizing operator that includes two additional parameters $a, b$ along with the regularization parameter $r$ \cite{kr1}:

\begin{equation}
\label{7}
f_r(t) = \int_0^\infty g(u) \Pi(r,a,b; tu) \mathrm{d} u,
\end{equation}

where $f_r(t)$ represents the regularized solution, and $f_r(t) \to f(t)$ as $r \to \infty$. The exact formula for the kernel $\Pi$ can be found in the referenced article or in software available on GitHub \cite{github}.

The accuracy of the computed regularized inverse $f_r(t)$ significantly depends not only on the optimal value of the regularization parameter but also on the appropriate values of additional method parameters. The author has proposed a heuristic criterion for determining acceptable values of all method parameters by minimizing the difference between two closely related regularized solutions $f_r^{(1)}$ and $f_r^{(2)}$ \cite{kr2}:

\begin{equation}
\label{8}
\min_{r, a, \alpha} \sum_{i=1}^{n}\left \{f_r^{(1)}(a, b, r; t_i) - f_r^{(2)}(a, b, r; t_i)\right \}^2.
\end{equation}

In cases where the additional parameters $a$ and $b$ are fixed, the comparison of two regularized solutions (\ref{8}) becomes a criterion for finding the regularization parameter.

\section{Parametrized regularization}

Calculation of the regularized inverse $f_r(t)$ with the help of formulae (\ref{7}), (\ref{8}) suggests that using semilogarithmic coordinates is more suitable for discretizing the integral (\ref{1})  \cite{github}. In this case, we have three discretization parameters $n, t_{\mathrm{min}}, t_{\mathrm{max}}$, representing the number of points $n$ logarithmically distributed on an integration interval $(t_{\mathrm{min}}, t_{\mathrm{max}})$.
As a result, the dimension and elements of the kernel matrix $K$ in equation (\ref{3}) will depend on discretization parameters.

Consequently, for computing the regularized solution of equation (\ref{3}) using Tikhonov regularization with a stabilizing matrix $\Omega$:

\begin{equation}
\label{81}
\min \left\{ ||K f -  g||^2 + \alpha ||\Omega f||^2 \right\},
\end{equation}

we need to compute the optimal values of the regularization parameter $\alpha$ as well as the optimal values of discretization parameters $n, t_{\mathrm{min}}, t_{\mathrm{max}}$.

From results described in the previous section, we can conjecture that the quasi-optimal values of all method parameters can be found by solving a minimization problem formulated based on a criterion for finding the regularization parameter. Then, for the chosen parametrized discretization, we encounter a four-dimensional minimization problem involving the regularization and discretization parameters.

In particular, the generalized cross-validation criterion \cite{Golub} is convenient in the case of a variable number of integration points. In this case, the quasi-optimal values of discretization and regularization parameters are obtained by solving the following minimization problem \cite{Hansen2}:

\begin{equation}
\label{9}
\min_{\alpha, n, t_{\mathrm{min}}, t_{\mathrm{max}}} {\frac{||K f_r -  g||^2}{\mathrm{trace}(I - KK^{\sharp})^2}},
\end{equation}

where the regularized solution $f_r = K^{\sharp} g$, and $I$ is the identity matrix.

The minimization problem for finding the quasi-optimal values of discretization and regularization parameters can also be formulated based on other criteria for searching the regularization parameter \cite{Tikhonov, Hansen2} and/or other approaches for parameterizing discretization.

As we will see in the next section, the quasi-optimal values of discretization and regularization parameters computed with the help of (\ref{9}) allow us to obtain quite satisfactory results and extract significantly more information from data in hand than conventional regularization techniques.

\section{Numerical Examples}

To illustrate the effectiveness of the proposed approach, we present a few examples of inverting Laplace transformations (\ref{1}) and solving problems of nuclear magnetic resonance relaxometry (\ref{100}). 

All examples were computed by using Tikhonov regularization (\ref{81}) with stabilizing matrices  $\Omega = I$ and  $\Omega = L_2$ , where $I$ is the identity matrix and $L_2$ is a matrix that approximates the second derivative of the solution. It is evident that matrix $L_2$ also depends on discretization parameters.

\begin{figure}[ht]
\begin{center}
\includegraphics[height=9.0 cm]{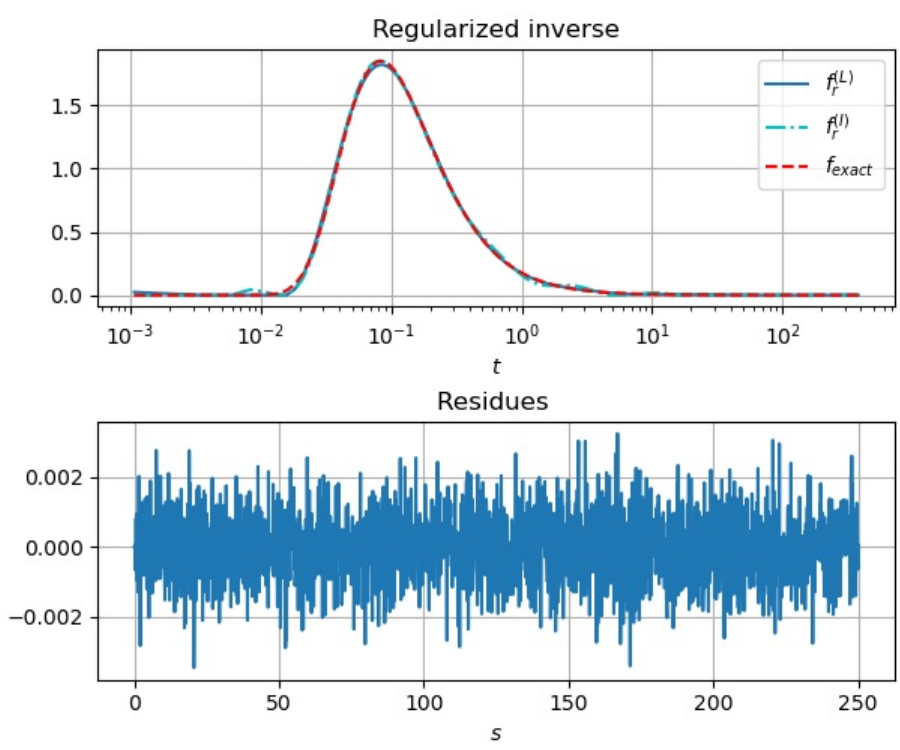}
\end{center}
\caption{\label{fig1} Inversion of a Laplace transform $g(s) = \exp(-\sqrt{s/2})$, $\sigma=0.001$.} 
\end{figure}

As a first example, let's consider the restoration of the pre-image function $ f(t) = \frac{1}{2 t \sqrt{2 \pi t} }\exp(-1/8t)$ from its Laplace transform $g(s) = \exp(-\sqrt{s/2})$ contaminated by adding normally distributed noise with $\sigma = 0.001$. As can be seen from the graph in Figure \ref{fig1}, we obtain a highly satisfactory solution for the selected level of noise.

The rest of the examples were constructed using a sum of log-normal distributions:

\begin{equation}
\label{12}
f(t) =  \sum_{i=1}^n a_i f(t;S_i, \theta_i),
\end{equation}

where 

\begin{equation}
\label{13}
f(t; S, \theta) = \frac{1}{t \sqrt{2\pi S }}\exp \left(\frac{\log{t}- \log{\theta}}{2 S}\right),
\end{equation}

For $S > 1$, the log-normal distribution reaches its maximum at $t \approx \theta$, and as $S \to \infty$, $f(t; S, \theta)$ tends to $\delta(t- \theta)$.

The image  $g(s)$ was computed numerically, and the experimental data have been simulated by adding normally distributed noise with  standard deviation $\sigma = const$.

\begin{figure}[ht]
\begin{center}
\includegraphics[height=9.0 cm]{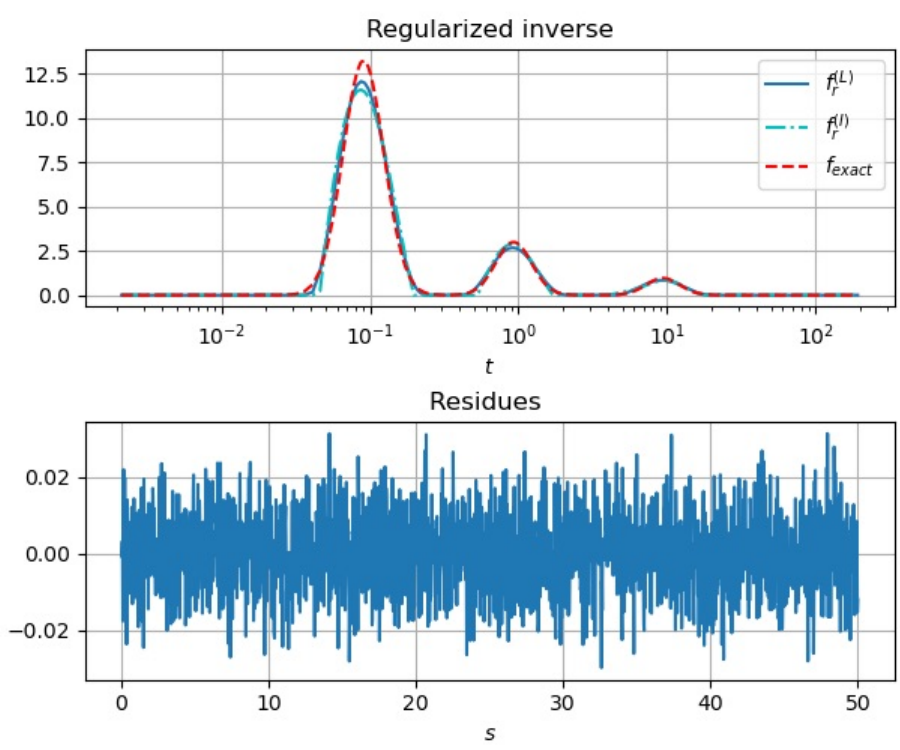}
\end{center}
\caption{\label{fig2} Restoration of a three-peak distribution. $\sigma = 0.01$, $a = [1, 2, 6], \theta = [0.1, 1, 10], S = [10, 13, 15]$ by inverse Laplace transformation}
\end{figure}

Graph in Figure \ref{fig2} shows restoration of a three-peak function from its Laplace transform. As can be seen, all peaks were resolved, and smoother peaks are restored more accurately. As observed, two regularized solutions are close to each other.

\begin{figure}[ht]
\begin{center}
\includegraphics[height=9.0 cm]{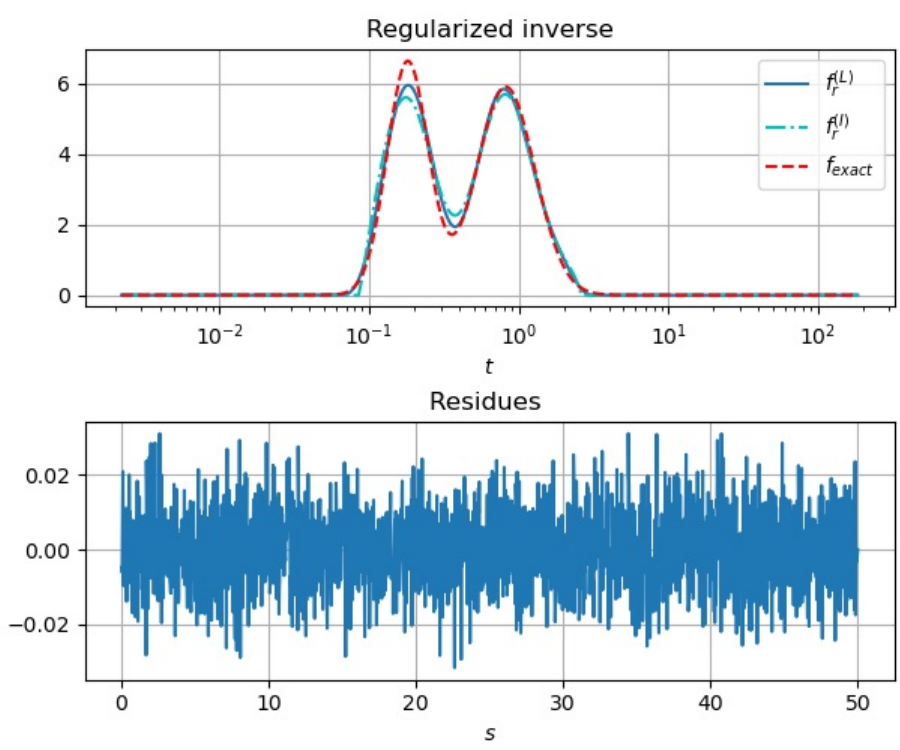}
\end{center}
\caption{\label{fig3} Restoration of a two-peak distribution. $\sigma = 0.01,  a = [1,  6], \theta = [0.2, 1],  S = [10, 5]$  by inverse Laplace transformation} 
\end{figure}

Graph in Figure \ref{fig3} shows the restoration of a sum of two slightly overlapping log-normal distributions. The results shown in the figure are quite satisfactory.

\begin{figure}[ht]
\begin{center}
\includegraphics[height=9.0 cm]{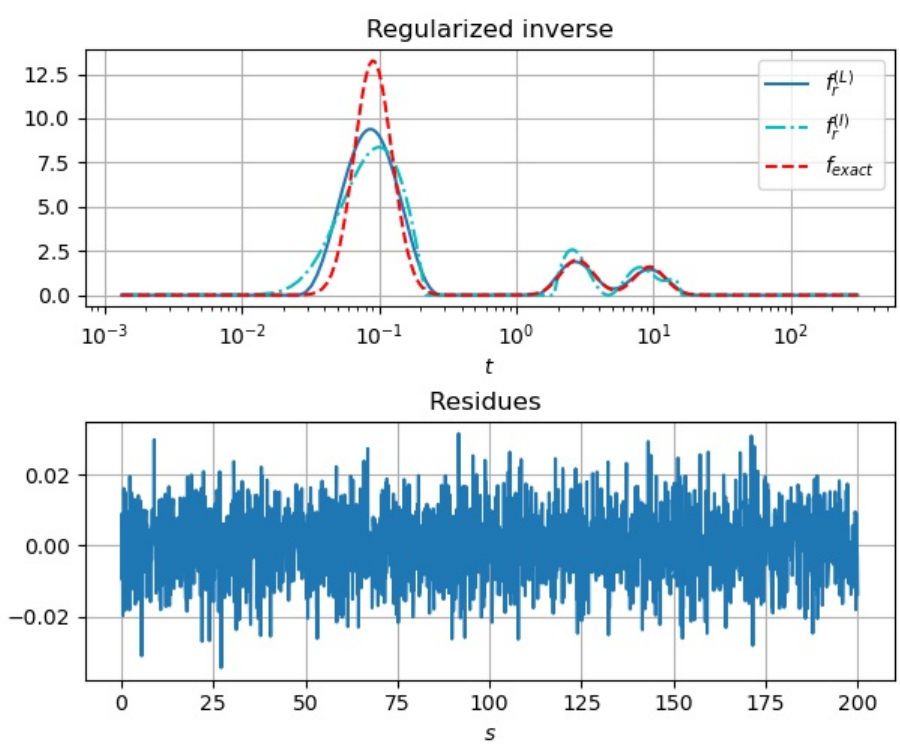}
\end{center}
\caption{\label{fig4} Deconvolution of NMR relaxation data. $\sigma = 0.01,  a = [1, 4, 10 ], \theta = [0.1, 3, 10],  S = [10, 13, 15]$} 
\end{figure}

\begin{figure}[ht]
\begin{center}
\includegraphics[height=9.0 cm]{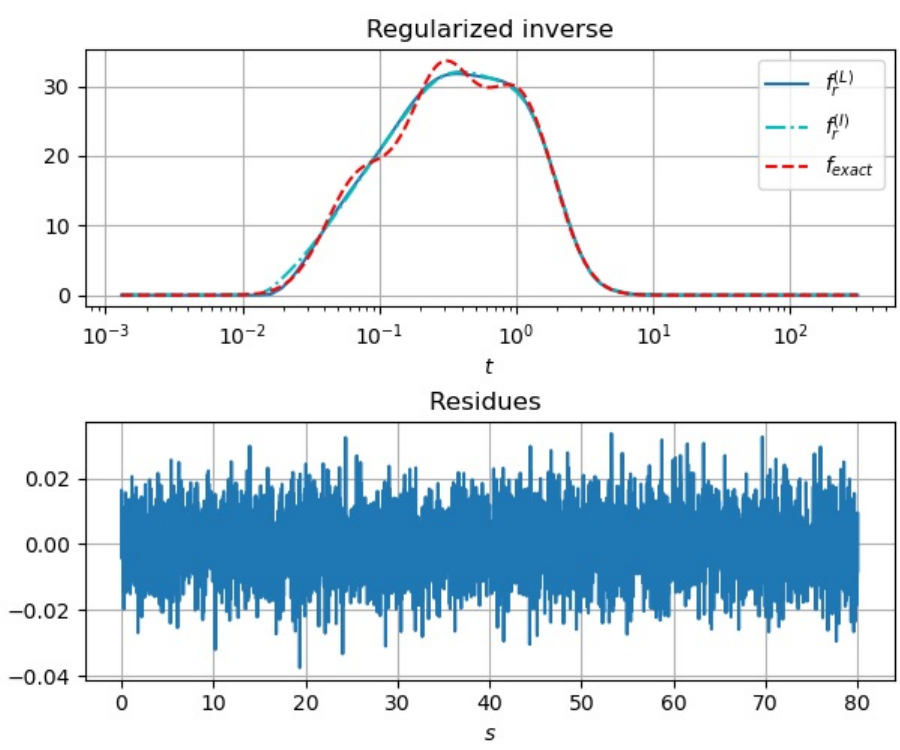}
\end{center}
\caption{\label{fig5} Deconvolution of NMR relaxation data. $\sigma = 0.01,  a = [2, 15, 50 ], \theta = [0.1, 0.4, 1.5],  S = [3, 3, 3]$} 
\end{figure}

Moving on to the application of the method in nuclear magnetic resonance relaxometry, Figure \ref{fig4} demonstrates the deconvolution of NMR relaxation data. We simulate a three-peak distribution as a sum (\ref{12}) with parameters $a = [1, 4, 10]$, $S = [10, 13, 15]$, $\theta = [0.1, 3, 10]$, and $\sigma = 0.01$. The results of deconvolving equation (\ref{100}) are shown in Figure \ref{fig4}, indicating successful restoration of all peaks. As can be seen, the kernel of the integral equation (\ref{100}) suppresses information about the function $f(t)$ for small $t$, consequently, the results for small $t$ are less reliable.  

The graph in Figure \ref{fig5} shows the results of deconvolution of NMR relaxation data in the case where components of a sum of distributions are highly overlapping.
As can be seen from the figure, in this case, the restored distribution smoothes out tiny details of the exact distribution. As observed, the closeness of solutions 
$f_r^{(I)}$ and $f_r^{(L)}$ does not necessarily mean the closeness to the exact solution.

\begin{figure}[ht]
\begin{center}
\includegraphics[height=9.0 cm]{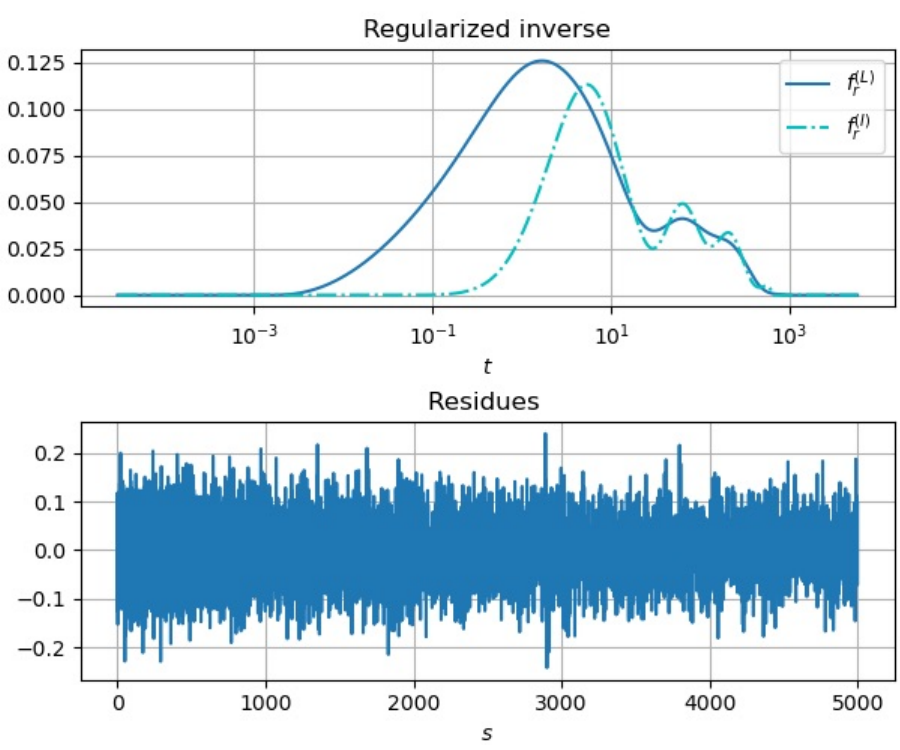}
\end{center}
\caption{\label{fig6} Deconvolution of measured NMR relaxation data.}
\end{figure}

Finally, figure \ref{fig6} showcases the deconvolution of experimentally measured NMR relaxation data. The noise-like residues depicted in the lower frame of the figure validate that the software performs as expected, providing a good fit to the data. The disagreement between two regularized solutions $f_r^{(I)}, f_r^{(L)}$ for small $t$ indicates that the results are not reliable in that area.

Additional testing results can be found in Jupyter notebooks available on Google Drive \cite{gd}.

\section{Conclusion}

In this paper, through examples of inversion of noisy Laplace transforms and NMR relaxation data, we have demonstrated that it is essential to introduce an appropriate parametrized discretization into the process of constructing regularizing operators for solving severely ill-posed problems. The optimal values for the regularization and discretization parameters can be computed simultaneously by solving a minimization problem formulated based on a regularization parameter search criterion.

Hence, it is reasonable to expect that other Fredholm integral equations of the first kind can be resolved more accurately by introducing an appropriate parametrized discretization.

\textbf{Acknowledgments}

The author is grateful to his wife, Helen Kryzhnyaya, for supporting his decades-long voluntary research on the inversion of real-valued Laplace transforms.

The author is also thankful to Green Imaging Technology, Inc. for providing experimental data for testing purposes and granting permission to use it in this paper.

\end{document}